\documentclass[12pt,oneside]{article}
\usepackage{harvard}
\usepackage{graphics,amsmath,amssymb}
\newcommand{\field}[1]{\mathbb{#1}}
\newtheorem{definition}{Definition}

\newtheorem{corollary}{Corollary}
\newtheorem{theorem}{Theorem}

\begin{document}
\begin{center}
\Large{On the history of the isomorphism problem of dynamical systems with special regard to \\
von Neumann's contribution}
\begin{verbatim}
\end{verbatim}
\large{Mikl\'{o}s R\'{e}dei (m.redei@lse.ac.uk) \\
and \\ Charlotte Werndl (c.s.werndl@lse.ac.uk)\footnote{This work is fully collaborative. Authors are listed alphabetically.}}
\begin{verbatim}
\end{verbatim}
Department of Philosophy, Logic and Scientific Method\\
London School of Economics\\
\begin{verbatim}
\end{verbatim}
\large{\emph{Forthcoming in: Archive for History of Exact Sciences} \\The final publication will be available at:
www.springer.com/mathematics/history+of+mathematics/journal/407
}
\begin{verbatim}
\end{verbatim}
\end{center}
\tableofcontents
\newpage

\section{Introduction}
The aim of this paper is to review some major episodes in the history of the spacial isomorphism problem of dynamical systems.  In particular, by publishing and analyzing (both systematically and in historical context) a hitherto unpublished letter\footnote{In Section \ref{Text} we reproduce von Neumann's letter in full. Doing so, we follow von Neumann's notation and terminology closely, and we keep his notation and terminology throughout the paper when we review subsequent developments about the isomorphism problem although some of that notation and terminology has become outdated. One instance of such non-standard usage of words is `spacial isomorphism' instead of `spatial isomorphism'.} written in 1941 by John von Neumann to Stanislaw Ulam, this paper aims to clarify von Neumann's contribution to discovering the relationship between spacial isomorphism and spectral isomorphism. The main message of the paper is that von Neumann's argument described (in a  sketchy way) in his letter to Ulam is the very first proof that spectral isomorphism and spacial isomorphism are \emph{not} equivalent because spectral isomorphism is weaker than spacial isomorphism: von Neumann shows that spectrally isomorphic  ergodic dynamical systems with mixed spectra need not be spacially isomorphic.

Dynamical systems are mathematical models of systems whose time evolution is deterministic, and these include physical systems described by classical mechanics and classical statistical mechanics. The origins of dynamical systems, especially the beginnings of ergodic theory, can be traced back to Boltzmann's work on the foundations of classical statistical mechanics in the second half of the 19th century \cite{Plato1991,Szasz1996}.

Modern ergodic theory developed in the early 1930s; the papers von Neumann \citeyear{Neumann1932ergodic2,vonNeumann1932,Neumann1932} and \citeasnoun{KoopmanNeumann1932} contributed substantially to the birth of this new discipline. More specifically, the isomorphism problem -- the problem of classifying the equivalence classes of dynamical systems with respect to a natural notion of spacial isomorphism -- was first formulated in \citeasnoun{Neumann1932}.

The development that gave a decisive boost to  dynamical systems theory, and which motivated in particular von Neumann to turn his attention to the study of the ergodic properties of dynamical systems, was the emergence of what has become called `the Koopman formalism':  dynamical systems theory in terms of functional analysis.  \citeasnoun{Koopman1931} observed that to each dynamical system there corresponds a group of unitary operators on the associated $L_2$ Hilbert space. Hence the properties of dynamical systems can be analyzed in terms of Hilbert space operator theory. The Koopman formalism proved to be very powerful and successful. Many stochastic properties of dynamical system (such as ergodicity and different types of mixing) could be fully characterized in terms of the Koopman operators associated with the dynamical system.  This success raised the hope that the isomorphism problem  could also be solved with this technique. Von Neumann \citeyear{Neumann1932} introduced a very natural notion of spectral isomorphism between Koopman representatives of dynamical systems, and he  proved that for a restricted class of ergodic dynamical systems (those with  pure point spectrum) spectral isomorphism is equivalent to spacial isomorphism. He also conjectured that this is true in  general. If this were true, it would have shown that the technique of the Koopman formalism is exhaustive in the sense of being strong enough to describe every  probabilistic aspect of dynamical systems.

It turned out,  however, that the situation is much more complicated, and a long and intricate history of research on the relation of spectral isomorphism to spacial isomorphism ensued after von Neumann's 1932 paper. A key event in this history is von Neumann's proof in his 1941 letter to Ulam that spacial isomorphism is stronger than spectral isomorphism. Von Neumann never published his proof but his result was referred to in important later publications on the problem. Given the significance of the isomorphism problem, von Neumann's letter is an important historical document.  Also, understanding von Neumann's proof is helpful in putting later developments into perspective. As it will be seen, von Neumann's presentation of the proof in his letter is sketchy, and we will fill in the details in order to make the proof understandable.

By giving a review of some major episodes in the history of the isomorphism problem we also hope to contribute to the extremely meager literature on the history of the isomorphism problem: while the literature on the history of ergodic theory is rich and there are a number of papers discussing specifically von Neumann's contribution to it (see \citename{Mackey1990} \citeyear*{Mackey1990}; \citename{Ornstein1990} \citeyear*{Ornstein1990}; \citename{Zund2002} \citeyear*{Zund2002} and the references therein), we are not aware of any  historically oriented paper devoted exclusively to the isomorphism problem, although some papers on the history of ergodic theory mention it (e.g.,  \citename{Mackey1990} \citeyear*{Mackey1990}; \citename{LoBello} \citeyear*{LoBello}).

The structure of the paper is the following. Section \ref{SIP} recalls some definitions, including the Koopman formalism, defines the spacial isomorphism problem and states the main  result on the problem that had been known before 1941 -- von Neumann's 1932 result. In Section \ref{Text} von Neumann's letter is reproduced without any comments on its content. Section \ref{Explanation} gives a detailed reconstruction of von Neumann's proof, filling in the (sometimes very large) gaps. Section \ref{Literature} reviews how von Neumann's result was referred to in the mathematics literature.  Section \ref{Kolm} is a collection of  important more recent results on the isomorphism problem. The paper is concluded with Section \ref{Conclusion} where we highlight the historical significance of von Neumann's letter.

\section{The Isomorphism Problem and von Neumann's 1932 Result }\label{SIP}

First, we recall some definitions that will be needed for the ensuing discussion, using  notation that is compatible with that in von Neumann's letter. $(\phi,\Sigma_{\phi},\mu,S)$ is a (discrete)\footnote{
We only focus on discrete-time systems because von Neumann's letter and all the relevant results are about them. (The only exception are the notions and results in \citeasnoun{Neumann1932} discussed in this and the next section. They are formulated for continuous-time systems but can easily be transferred to discrete-time systems.)} \textit{dynamical system} if $\phi$\footnote{`$\phi$' is the notation von Neumann uses in his letter for the phase space (cf.\ Section~\ref{Text}). This differs from modern usage (where `$\phi$' is usually reserved for functions).}
is a set (called the phase space), $\Sigma_{\phi}$ is a $\sigma$-algebra of subsets of $\phi$, $\mu$ is a measure on $\Sigma_{\phi}$ with $\mu(\phi)=1$, and $S:\phi\rightarrow\phi$  is a bijective measurable function such that also $S^{-1}$ is measurable.\footnote{We also assume, as is standard in ergodic theory, that $(\phi,\Sigma_{\phi},\mu)$ is a Lebesgue space (cf.\ \citename{Petersen1983} \citeyear*{Petersen1983}).
$(\phi,\Sigma_{\phi},\mu)$ is called a \emph{Lebesgue space} if either $\phi$ is countable or there is a measure space $(K,\Sigma_{K},\nu)$, where $K=[a,b)\subseteq \field{R}$ is an interval, $\Sigma_{K}$ is the Lebesgue $\sigma$-algebra and $\nu$ the Lebesgue measure, such that the following holds:
There is a countable set $W\subseteq \phi$, there is a $K_{0}\subseteq K$ with $\nu(K_{0})=1$, there is a $\phi_{0}\subseteq \phi$ with $\mu(\phi_{0})=1$ and there is a bijective function $c:\phi_{0}\setminus{W}\rightarrow K_{0}$ such that \begin{description} \item[(i)] $c(A)\!\in\!\Sigma_{K}$ for all $A\!\in\!\Sigma_{\phi}$, $A\subseteq\phi_{0}\setminus W$, \item[(ii)] $c^{-1}(B)\in\Sigma_{\phi}$ for all $B\in\Sigma_{K},B\subseteq K_{0}$, \item[(iii)] $\nu(c(A))=\mu(A)$ for all $A\in\Sigma_{\phi},\,A\subseteq\phi_{0}\setminus W$.
\end{description}
}

The definition of spacial isomorphism captures the idea that dynamical systems are equivalent from a probabilistic perspective.
\begin{definition}\label{def:spacial_isomorphism}
{\rm
Two dynamical systems $(\phi,\Sigma_{\phi},\mu,S)$ and
$(\psi,\Sigma_{\psi},\nu,T)$ are called \textit{spacially isomorphic} if there are measurable sets $\phi_{0}\subseteq \phi$ and $\psi_{0}\subseteq \psi$ with the properties 1. and 2. below and a bijection $c:\phi_{0}\!\rightarrow\!\psi_{0}$ satisfying 3.{\bf (i)}-{\bf(iii)}:
\begin{enumerate}
\item $S(\phi_{0})\subseteq\phi_{0}$, \quad $T(\psi_{0})\subseteq\psi_{0}$;
\item $\mu(\phi\setminus\phi_{0})=0$, \quad $\nu(\psi\setminus \psi_{0})=0$;

\item \begin{description}
  \item[(i)] $c$ and $c^{-1}$ are measurable with respect to \\
  $\Sigma_{\phi_{0}}:=\{A\cap\phi_{0}\,|\,A\in\Sigma_{\phi}\}$ and
$\Sigma_{\psi_{0}}:=\{B\cap\psi_{0}\,\,|\,\,B\in\Sigma_{\psi}\}$,
  \item[(ii)] $\nu(c(A))=\mu(A)$ for all
$A\in\Sigma_{\phi_{0}}$,
  \item[(iii)] $c(S(x))=T(c(x))$ for all $x\in\phi_{0}$.
\end{description}
\end{enumerate}
}
\end{definition}

Spacial isomorphism was first  defined in \citeasnoun{Neumann1932}. Von Neumann's \citeyear{Neumann1932} paper also was the first that called for a \emph{classification of dynamical systems up to spacial isomorphism} -- this is the \emph{spacial isomorphism problem}. The spacial isomorphism problem is widely regarded as one of the most important problems, if not \emph{the} most important problem, of ergodic theory (cf.\ \citename{Halmos1956} 1956, 96; \citename{Halmos1961} 1960, 75; Petersen 1983, 4; \citename{Rohlin1960} 1960, 1). This is understandable because ergodic theory is the theory of the probabilistic  behaviour of dynamical systems; it is thus a crucial question which systems are equivalent from the probabilistic perspective.

There are two main ways to approach this problem. First, one tries to find \emph{spacial invariants} (properties which are the same for spacially isomorphic dynamical systems) which show that certain kinds of dynamical systems are not spacially isomorphic. Second, one aims to find a sufficiently large number of invariants which, taken together, provide sufficient conditions for systems to be isomorphic.

Next we briefly describe the Koopman formalism, which will lead to the notion of spectral isomorphism of dynamical systems.
Given a dynamical system $(\phi,\Sigma_{\phi},\mu,S)$ consider the Hilbert space $L_{2}(\phi)$ of complex-valued square integrable functions on $\phi$. The operator
$$
U:L_{2}(\phi)\rightarrow L_{2}(\phi), \qquad \,U(f)=f(S(x))$$
is called the \textit{Koopman operator} of $(\phi,\Sigma_{\phi},\mu,S)$. \citeasnoun{Koopman1931} showed that the Koopman operator is \textit{unitary} \citeaffixed{Arnold1968}{cf.}. The spectrum of the Koopman operator $U$ is called the \emph{spectrum} of the dynamical system.

\begin{definition}\label{def:spectral_isomorphism}
{\rm
Two dynamical systems $(\phi,\Sigma_{\phi},\mu,S)$ and $(\psi,\Sigma_{\psi},\nu,T)$
are  called \textit{spectrally (or unitary) isomorphic} if their Koopman operators $U$ and $V$ are unitarily equivalent, i.e.\ if there is a unitary operator $W$ between $L_2(\phi)$ and $L_2(\psi)$ such that
$U=W^{*}VW$.
}
\end{definition}
\emph{Spectral invariants} are properties which are the same for spectrally isomorphic dynamical systems. Obviously, the spectrum of a dynamical system is a spectral invariant: two dynamical systems cannot be spectrally isomorphic if they have different spectrum.  However, having the same spectrum does not entail that they are spectrally isomorphic (\citename{Halmos1951} 1951, 75).

\emph{What is the relation of spectral and spacial invariants? Are they identical? }If yes, then spacial isomorphism is equivalent to spectral isomorphism. These questions were raised by von Neumann \citeyear{Neumann1932}. An affirmative answer to the question about the equivalence of spectral and spacial isomorphism would mean that  the Koopman formalism can deal with all the probabilistic questions about dynamical systems (cf.\ \citename{Neumann1932} \citeyear*{Neumann1932}; \citename{Halmos1957} 1957).

The main result on this problem before von Neumann's 1941 letter to Ulam was presented in von Neumann's \citeyear{Neumann1932} seminal paper. The paper starts by listing several properties of dynamical systems which are spacial invariants but are also spectral invariants and hence can be characterized completely in terms of the Koopman representatives of dynamical systems. We recall here two of the most prominent of these properties and their spectral characterizations.

A dynamical system $(\phi,\Sigma_{\phi},\mu,S)$ is called
\citeaffixed{Petersen1983}{cf.}
\begin{description}
\item[Ergodic:] (or \textit{metrically transitive}): if it does not have non-trivial $S$-invariant measurable sets, i.e.\ if there is no set $A\in\Sigma_{\phi}$, $0<\mu(A)<1$, such that $S(A)=A$.\\
    \emph{Spectral characterisation} of ergodicity: $1$ is a simple proper value\footnote{`Proper value' and `proper function' are also called `eigenvalue' and `eigenfunction', but we follow here (and throughout the paper) von Neumann's terminology in his letter (see Section \ref{Text}). `Proper value' and `proper function' also is the terminology in the literature we will refer to (e.g., Halmos \citeyear*{Halmos1949}, \citeyear*{Halmos1956}). A proper value $\delta$ is simple (also called `non-degenerate') if the proper functions belonging to $\delta$ span a one-dimensional subspace in $L_2(\phi)$.\label{footnote:propervalue}} of the Koopman operator $U$.
\item[Weakly mixing:] if for all $A,B\in\Sigma_{\phi}$ we have
$$
\lim_{t\rightarrow\infty}\frac{1}{t}\sum_{i=0}^{t-1}|\mu(S^{i}(A)\cap B)-\mu(A)\mu(B)|=0.
$$
\emph{Spectral characterisation} of weak mixing: 1 is the only proper value of the Koopman operator and this proper value is simple.
\end{description}

The focus of von Neumann \citeyear*{Neumann1932} then turned to the question about the relationship between spacial and spectral isomorphism. It is obvious that if $(\phi,\Sigma_{\phi},\mu,S)$ and
$(\psi,\Sigma_{\psi},\nu,T)$ are spacially isomorphic, then they are spectrally isomorphic: if $c$ is the spacial isomorphism, then $W$ defined by
$$
W(f(\phi)) = f(c^{-1}(\psi))
$$
sets up the spectral equivalence of the dynamical systems. This implies that spectral invariants, such as the spectrum of a system, are also spacial invariants. Thus the question to answer was whether the converse is true.

To tackle this question, \citeasnoun{Neumann1932} distinguished between several types of spectra of dynamical systems. A dynamical system $(\phi,\Sigma_{\phi},\mu,S)$ is said to have \textit{pure point spectrum} if the proper functions of the Koopman operator $U$ form a basis of $L_{2}(\phi)$. A dynamical system has \emph{pure continuous spectrum} if the only proper functions of its Koopman operator are $\delta 1$, $\delta\in\field{C}$. A dynamical system has \emph{mixed spectrum} if it has neither a pure point spectrum nor a pure continuous spectrum.

Von Neumann then proved the following theorem.\footnote{For a modern proof, see \citename{Halmos1957} (1957, 46--50). This result is also discussed by the more historically oriented papers \citeasnoun{Halmos1949}, \citeasnoun{Halmos1957}, \citename{Mackey1974} (1974, 197), \citename{Mackey1990}~\citeyear{Mackey1990}, \citename{ReedSimon1980} (1980, Section VII.4),
\citename{Rohlin1967} (1967, 2 and 5) and \citename{Weiss1972} (1972, 672--673).}
\begin{theorem}\label{T1}
Let $(\phi,\Sigma_{\phi},\mu,S)$ and $(\psi,\Sigma_{\psi},\nu,T)$ be ergodic dynamical systems with pure point spectrum. Then the dynamical systems are spacially isomorphic if, and only if, they are spectrally isomorphic \cite{Neumann1932}.
\end{theorem}
It is easy to see that spectrally isomorphic systems have the same proper values. Thus Theorem~\ref{T1} tells us that ergodic systems with pure point spectrum are spacially isomorphic precisely when they have the same proper values.

Theorem~\ref{T1} was the main result on the isomorphism problem before von Neumann' 1941 letter. The theorem shows that for ergodic systems with pure point spectrum the Koopman formalism can deal with all  probabilistic questions about dynamical systems. To the best of our knowledge, nothing was known about the other cases of a mixed spectrum and continuous spectrum before the 1941 letter. Yet in these early days of ergodic theory it was sometimes conjectured that, in general, dynamical systems are spacially isomorphic just in case they are spectrally isomorphic (cf.\ \citename{Halmos1961} \citeyear*{Halmos1961}, 77). Even \citeasnoun{vonNeumann1932}, when discussing systems with pure continuous spectrum, states that results similar to Theorem~1 can be expected in this case.

Let us now turn to the next step in the history of the isomorphism problem, viz.\ von Neumann's letter.

\section{Von Neumann's 1941 Letter: the Text}\label{Text}

\begin{figure}
\centering
\includegraphics{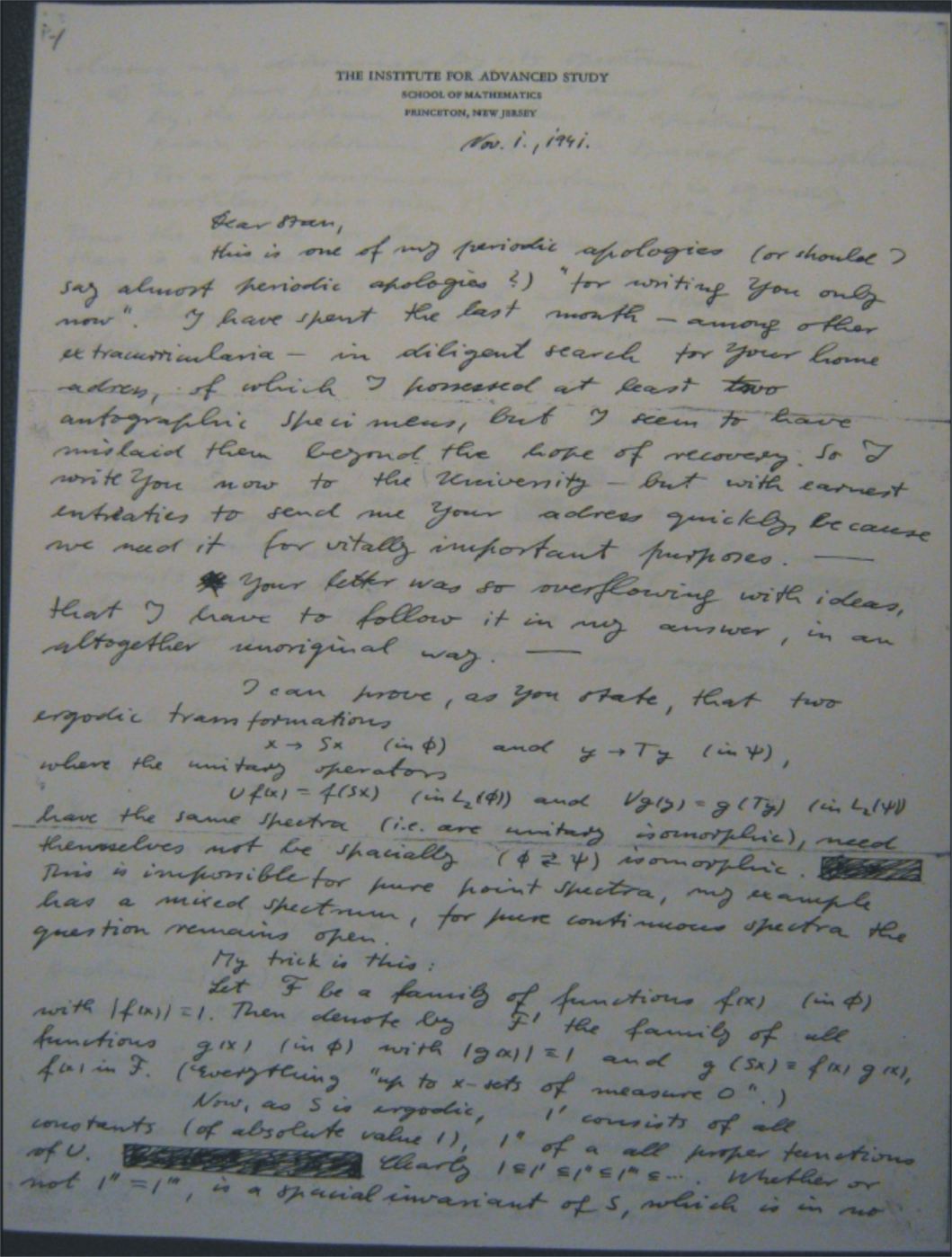}
\caption{\normalsize{First page of von Neumann's 1941 letter}}
\end{figure}
The letter below was written by John von Neumann to Stanislaw Ulam (the first page of the original is shown in Figure~1). Von Neumann and Ulam were not only colleagues but friends as well. This explains the informal tone and chatty style of the letter and the references to personal matters. Ulam and von Neumann corresponded on a number of issues (both scientific and non-scientific) over many years. A selection of von Neumann's letters to Ulam was published in \citeasnoun{Neumann2005}. One of these, not fully dated but most likely written in 1939, touches upon the isomorphism problem briefly \cite[252--253]{Neumann2005}:  Ulam's attention is called to the paper \citeasnoun{Neumann1932}, and in particular to the theorem that systems with pure point spectrum are spacially isomorphic if, and only if, they are spectrally isomorphic (Theorem 1).
This letter suggests that von Neumann did not yet know in 1939 that spacial and spectral isomorphism do not always go together and, in particular, that he had not yet discovered the result described in the 1941 letter.

To the best of our knowledge, the 1941 letter has not been published. The original, hand-written version is in the S.~Ulam papers, in the Archives of the American Philosophical Society, Philadelphia, U.S.A. The language of the letter is English. We reproduce it below without comments and keeping the original spelling. The content of it will be analyzed in detail in the next section.

\begin{center}

-------------------------------

\end{center}

\bigskip
November 1, 1941.\\

Dear Stan,

\noindent this is one of my periodic apologies (or should I say almost periodic apologies?) ``for writing you only now''. I have spent the last month -- among other extracurricularia -- in diligent search for your home adress, of which I possessed at least two autographic specimens, but I seem to have mislaid them beyond the hope of recovery. So I write you now to the University -- but with earnest entreaties to send me your adress quickly, because we need it for vitally important purposes. ------

Your letter was so overflowing with ideas, that I have to follow it in my answer, in an altogether unoriginal way. ------

I can prove, as you state, that two ergodic transformations
\begin{equation}
x\rightarrow Sx\,\,\,(\textnormal{in}\,\,\phi)\,\,\,\textnormal{and}\,\,\,y\rightarrow Ty\,\,\,(\textnormal{in}\,\,\psi),
\end{equation}
where the unitary operators
\begin{equation}
Uf(x)=f(Sx)\,\,\,(\textnormal{in}\,\,L_{2}(\phi))\,\,\,\textnormal{and}\,\,\,Vg(y)=g(Ty)\,\,\,(\textnormal{in}\,\,L_{2}(\psi))
\end{equation}
have the same spectra (i.e.\ are spectrally isomorphic), need themselves not be spacially $(\phi \leftrightarrow \psi)$ isomorphic. This is impossible for pure point spectra, my example has a mixed spectrum, for pure continuous spectra the question remains open.

My trick is this:

Let $F$ be a family of functions $f(x)$ (in $\phi$) with $|f(x)|=1$. Then denote by $F'$ the family of all functions $g(x)$ (in $\phi$) with $|g(x)|=1$ and $g(Sx)=f(x)g(x)$, $f(x)$ in $F$. (Everything ``up to $x$-sets of measure $0$''.)

Now,\footnote{In this paragraph, and throughout von Neumann's letter, 1 stands for the constant function on $\phi$, taking the value $1\in\mathbb{R}$. This is slightly confusing because the symbol refers both to a function and the real number `one'; however, the context will always make clear what 1 stands for. Also note another slight abuse of notation in this paragraph and in the letter: The previous paragraph defines $F'$ for a {\em set\/} $F$ of functions. However, von Neumann shortens $\{1\}'$ to $1'$ in case where $F=\{1\}$, and the same holds for higher primes. So, for instance, $\{\{1\}'\}'$ is denoted by $1''$ etc. \label{primeabuse}} as $S$ is ergodic, $1'$ consists of all constants (of absolute value 1), $1''$ of all proper functions of $U$. Clearly, $1\subseteq 1'\subseteq 1''\subseteq 1'''\subseteq\ldots$. Whether or not $1''=1'''$, is a spacial invariant of $S$, which is in no obvious way determined by its spectrum. But:

\begin{description}
\item
$\alpha$) For a pure point spectrum it must be determined by the spectrum, since then the spectrum is known to determine $S$ up to a spacial isomorphism.
\item $\beta$) For a pure continuous spectrum it is equally worthless, since then $1'=1''$, hence $1''=1'''$.\end{description}
Thus the chance for this invariant arises when there is a mixed spectrum.

Now let $\phi$ be the space of all $u,v$ (both mod$1$), with Lebesgue measure. Choose a fixed, irrational $\gamma$, and define
\begin{equation}
S(u,v)=(u+\gamma,v+u).
\end{equation}
(This example was constructed sometime ago by Halmos, for a different purpose.)

It is easy to show, that $S$ has this spectrum
\begin{description}
\item I) A simple point spectrum $e^{2\pi i k \gamma}$, $k=0,\pm 1,\pm 2,\ldots$.
\item II) An everywhere infinitely multiple, Lebesgue measure spectrum, covering all of $|\lambda|=1$
\end{description}
$1''$ consists precisely of all $e^{2\pi i k u}$, $k=0, \pm 1, \pm 2, \ldots$. $1'''$ contains $e^{2\pi i v}$. Hence $1''\neq 1'''$.

Consider on the other hand any ergodic transformation
\begin{equation}
w \rightarrow Rw\,\,\,(\textnormal{in}\,\,\Omega)
\end{equation}
which has this spectrum:
\begin{description}
\item I') A simple point spectrum $1$.
\item II') Same as II) above.
\end{description}
(E.g.\ the Hopfe one: $w=(x_{n}\,|\,n=0,\pm 1,\pm 2, \ldots)$, $Rw=(x_{n+1}\,|\,n=0,\pm 1,\pm 2, \ldots)$ --- each $x_{n}$-s range being $\pm 1$.)\\
Define
\begin{equation}
T(u,w)=(u+\gamma,Rw).\end{equation}

Then it is easy to show that $T$ has the same spectrum I)+II) as $S$. And a simple computation gives $I'''=I''=$ set of all $e^{2\pi i k u}$, $k=0,\pm 1,\pm 2,\ldots$. I.e., $1''=1'''$. Thus $S$ and $T$ are spacially non-isomorphic.

Re Kakutani: K. has decided to stay on. He has an I.A.S. stipend for 1941/42. If the war should force him to remain beyond Sept.\ 1942, we will do everything humanly possible for him. If he could be elected to the Harvard Society of Fellows, that would be super-excellent. Did you find out from Henderson whether there is any chance of this? It will be too glad to do anything, to write to anybody, etc. -- if there is a chance of achieving something. ------

Norberts resignation from the Academy: I have no idea whether he has resigned or not. I remember that he used to talk about it 2--3--4 years ago. I think it is nonsense. ------

Many thanks for the invitation to talk at Madison, including the Mammon. I have no way to foresee the future in this slightly opaque century -- but I hope that I may get Midwest again before too long -- but I don't know for the moment when or how. ------

Many congratulations to your + Oxtobys early paper, which appeared so late in the ``Annals''. It is really a pleasure to read it. ------

My respectful homage to Fran\c{c}oise, God knows that my similarity to Riquet \'{a} la Houppe (if any -- which I doubt) is rigorously exterior. I am a miserable sinner, and I never rescue maidens, except for base and egoistical motives. But her high -- and utterly undeserved -- consideration thrills me to the core of my being. ------

I am also highly pleased by the Russian campaign, although I fear that it will be like Old Times -- i.e.\ 1914--1917 -- when the German wore Russia down after 2--3 years. But they will probably not be able to get rid of some kind of war in that direction -- more or less like honorable Sino-Japanese incident.

I suppose that Schickelgruber is licked alright, but it will be a long and bloody affair. And as to the US: You know what the court physician said to the german prince when his wife bore him a daughter: ``Majest\"{a}t werden sich nochmal bem\"{u}hen m\"{u}ssen.'' ------

Klari hopes for an epistle from you. It seem to be a difficult piece of accounting-letters. She has had a rather unpleasant eye-trouble lately -- probably over-strain -- and it is still not quite over. ------

I am writing a book on games and economics with the economist O.\ Morgenstern. Whadayasay.

Also a lot whatnots on operator theory, for which I have a monopoly in production as well as in consumption. So that's that.

Please remind me to Fran\c{c}oise.

With the best from both of us

hoping to hear from you soon

as ever
John.

\begin{center}

-------------------------------

\end{center}

\section{Von Neumann's 1941 Letter: Explanation of the Proof}\label{Explanation}

\emph{In the letter von Neumann shows that spectrally isomorphic dynamical systems with a mixed spectrum need not be spacially isomorphic.} He starts with the following construction (see paragraphs 4--5). Consider an ergodic dynamical system $(\phi,\Sigma_{\phi},\mu,S)$. Recall that $L_{2}(\phi)$ is the set of all complex-valued square-integrable functions on $\phi$. Given a family $F$ of functions $f \in L_{2}(\phi)$ whose absolute value is always $1$ (i.e.\ $|f(x)|=1$ for all $x\in\phi$), let $F'$ be the family of functions $g \in L_{2}(\phi)$ whose absolute value is always $1$ (i.e.\ $|g(x)|=1$ for all $x\in\phi$) where there is an $f\in F$ such that
\begin{equation}\label{generalised}
g(Sx)=f(x)g(x).
\end{equation}
Let $1$ denote the function $f(x)=1$ for all $x\in\phi$, and consider $1'$ (see footnote \ref{primeabuse}). For ergodic dynamical systems the only \textit{invariant functions} (i.e.\ functions $g$ with $g(S(x))=g(x)$) are the constants \cite{Arnold1968}.  Hence $1'$ is the set of all constant functions of absolute value 1. Thus $1''$ is the set of all proper functions of $U$ of absolute value $1$. Clearly, for any sets $G,H$, if $G\subseteq H$, then $G'\subseteq H'$. Let us call the functions in $1''$, $1'''$, $1''''\ldots$ \textit{generalised proper functions} with \textit{generalised proper values} $1'$, $1''$, $1'''\ldots$.\footnote{In the literature they are called like this \cite{Abramov,Halmos1949,Halmos1956,Rohlin1960}. They have been first introduced when proving the result of von Neumann's letter (see Section~\ref{Literature}). This terminology is motivated as follows: $1''$ are proper functions of absolute value $1$ with proper values $1'$ of absolute value $1$, and the generalised proper functions and proper values are obtained by allowing $f$ in equation~(\ref{generalised}) to be not just the constant functions but any function in $1''$ or $1'''$ etc. From unitarity it follows that the absolute value of any proper value of the Koopman operator is 1. It suffices to consider proper functions of absolute value $1$ because for ergodic systems the absolute value of any proper function $g$ is constant, i.e.\ $|g(x)|=c$ for a constant $c$ for all $x\in\phi$ \cite[26]{Arnold1968}.}
Because $1\subseteq 1'$, $1\subseteq 1'\subseteq 1''\subseteq 1'''\ldots$.

Von Neumann states (see paragraphs 3 and 5) that whether $1''=1'''$ is in no obvious way determined by the spectral properties. What is clear is that for spacially isomorphic dynamical systems $(\phi,\Sigma_{\phi},\mu,S)$ and $(\psi,\Sigma_{\psi},\nu,T)$ either for both $1''=1'''$ or for both $1''\neq 1'''$. This is so because for spacially isomorphic systems their unitary operators are related via $U=W^{*}VW$ where $W(f)=f(c^{-1}(y))$. Therefore, $U(g)=fg$ if, and only if, $V(W(g))=W(fg)=W(f)W(g)$, and hence $g$ is a generalised proper function of $(\phi,\Sigma_{\phi},\mu,S)$ if, and only if, $W(g)$ is a generalised proper function of $(\psi,\Sigma_{\psi},\nu,T)$.

Von Neumann's idea now was \textit{to find two dynamical systems which are spectrally isomorphic but where for one system $1''=1'''$ and for the other $1''\neq 1'''$, implying that they are not spacially isomorphic.}
Note that this strategy will not work for ergodic systems with pure point spectrum because for those systems spacial and spectral isomorphism are equivalent.  Spectrally isomorphic systems have the same  proper values; hence this strategy will also {\em not} work for systems with {\em pure continuous spectrum} because their only proper functions are $\delta 1$, $\delta\in\field{C}$, and thus $1'=1''=1'''$. Hence von Neumann focused on systems with \textit{mixed spectrum} and constructed two spectrally isomorphic dynamical systems having mixed spectrum such that for one $1''=1'''$ holds, and for the other $1''\neq 1'''$ holds.

Consider (see paragraph 6) the dynamical system where $\phi=[0,1)\times [0,1)$, $\Sigma_{\phi}$ is the Lebesgue $\sigma$-algebra, $\mu$ is the Lebesgue measure and
\begin{equation}
S(u,v)=(u+\gamma\,(\textnormal{mod}\,1),v+u\,(\textnormal{mod}\,1))
\end{equation}
for a fixed irrational $\gamma$. It is shown for this system that $1''\neq 1'''$.

It is well known (see paragraph 7) that the functions
$$
g_{k,m}(u,v)=e^{2\pi i k u}e^{2\pi i v m} \qquad k,m\in\field{Z}
$$
form an orthonormal basis of $L_{2}(\phi)$ \cite{Arnold1968,Halmos1956}. We have
$$
U(g_{k,m})=e^{2\pi i k \gamma}g_{k+m,m};
$$
hence $g_{k,0}$, $k\in\field{Z}$, are proper functions. Let $H_{0}$ be the closed subspace in $L_{2}(\phi)$ generated by $\{g_{k,0}\ : \ k\in\field{Z}\}$. The functions $h_{k}:=g_{k,0}$, $k\in\field{Z}$, form an orthonormal basis of $H_{0}$. Thus $U$ restricted to $H_{0}$ has pure point spectrum with proper functions $h_{k}$ and proper values $e^{2\pi i k \gamma}$  -- this is point I)  referred to by von Neumann in the letter.

Now by applying $U$, the functions $g_{k,m}$ for $m\neq 0$ get permuted among themselves and multiplied by constant factors. We get rid of these constant factors by setting
\begin{equation} f_{k,m}=a_{k,m}g_{k,m},\,\,\textnormal{with}\,\, a_{k,m}=(\sqrt[8m]{e^{2\pi i \gamma k}})^{(2n-m)^{2}},\,\,k,m\in\field{Z},\,m\neq 0.
\end{equation}
A simple calculation yields $U(f_{k,m})=f_{k+m,m}$. Since
$|a_{k,m}|=1$, the functions $f_{k,m}$ are still an orthonormal set.
We can label the functions $f_{k,m}$ such that we obtain a sequence of functions $h_{i,j}$, $i,j\in\field{Z}$, with $U(h_{i,j})=h_{i+1,j}$. The subspace generated by $h_{i,j}$ is the orthogonal complement $H_0^{\perp}$ of $H_0$ in $L_{2}(\phi)$, and $h_{i,j}$ is an orthonormal basis  of $H_0^{\perp}$.

An operator $O$ is said to have an \textit{infinitely multiple homogeneous Lebesgue measure spectrum} if there exists an orthonormal basis $b_{i,j}$, $i,j\in\field{Z}$, with $O(b_{i,j})=b_{i+1,j}$
(\citename{Arnold1968}~\citeyear*{Arnold1968}, 28--30; \citename{Cornfeldetal1982}~\citeyear*{Cornfeldetal1982}, Appendix~2).
Consequently, $U$ acting on $H_0^{\perp}$ has an infinitely multiple homogeneous Lebesgue measure spectrum -- this is point II) referred to by von Neumann.

Let us show that there are no other proper functions of absolute value $1$ besides $\delta g_{k,0}$, $|\delta|=1$. Suppose that $g(u,v)$ is a proper function with proper value $\lambda$. Because $g_{k,m}$ is an orthonormal basis, we can expand $g$ as
$$
g(u,v)=\sum_{k,m}\chi_{k,m}g_{k,m}(u,v), \qquad \chi_{k,m}\in\field{C}.
$$
Then
\begin{equation}
Ug=\sum_{k,m}e^{2\pi i (k-m)\gamma}\chi_{k-m,m}g_{k,m}=\lambda\sum_{k,m}\chi_{k,m}g_{k,m}.
\end{equation}
Thus $|\chi_{k-m,m}|=|\lambda\chi_{k,m}|$ for all $k,m$.  It is well known (e.g., Halmos \citeyear*{Halmos1951}, 20) that if $\kappa,\kappa_{j}$ are elements in an Hilbert space with $\kappa_{j}$ pairwise orthogonal, then $\sum_{j}\kappa_{j}=\kappa$ if, and only if, $\sum_{j}\left\|\kappa_{j}\right\|^{2}$ converges. Hence  $\sum_{k,m}|\chi_{k,m}|^{2}$ converges, implying that $\chi_{k,m}=0$  whenever $m\neq 0$.
Therefore, $1''$ -- the set of proper functions -- consists of $\delta e^{2\pi i k u}$, $|\delta|=1$, $k\in\field{Z}$. Because for $g=e^{2\pi i v}$ we have
$$g(S(u,v))=e^{2\pi i(u+v)}=e^{2\pi i u}g;
$$
$1'''$ contains $e^{2\pi iv}$. Consequently, $1''\neq 1'''$.

Consider (see paragraph 8) a dynamical system $(\Omega,\Sigma_{\Omega},\rho,R)$ with Koopman operator $X$  where
\begin{description}
\item[(i)]
$\delta 1$, $|\delta|=1,\,\delta\in\field{C}$, are the only proper functions with absolute value $1$;
\item[(ii)] $X$ acting on $H_{1}^{\perp}$, where $H_{1}$ is the subspace generated by the function $1$, has infinitely multiple homogeneous Lebesgue measure spectrum.
\end{description}
For instance, let $\Omega$ be the set of all bi-infinite sequences $w=(...w_{-1},w_{0},w_{1}...)$ with $w_{i}$ either $-1$ or $1$. Let $\Sigma_{\Omega}$ be the $\sigma$-algebra
generated by the semi-algebra of sets
\begin{equation}
C^{e_{1}...e_{n}}_{i_{1}...i_{n}}\!\!=\!\!\{w\in \Omega\,|\,w_{i_{1}}\!\!=\!\!e_{1},...,w_{i_{n}}\!\!=\!\!e_{n},\,i_{j}\!\in\!\field{Z},\,i_{1}\!\!<...<\!\!i_{n},\,e_{j}\!\in\!\{-1,1\},\!\,1\!\leq j\!\leq n\}.
\end{equation}

Let $\rho$ be the unique extension to a measure on $\Omega$ of the pre-measure defined by
$$
\bar{\rho}(C^{e_{1}...e_{n}}_{i_{1}...i_{n}})=p_{e_{1}}...p_{e_{n}},
$$
where $p_{-1}$ and $p_{1}$ are real numbers with $p_{-1}+p_{1}=1$, $0<p_{-1},p_{1}<1$. Finally, let $R$ be defined by
$$
R(...w_{-1},w_{0},w_{1}...)=(...w_{0},w_{1},w_{2}...).
$$
The dynamical system $(\Omega,\Sigma_{\Omega},\rho,R)$ is called a \textit{Bernoulli shift} on the symbols $-1$ and $1$ \citeaffixed{Werndl2009c}{cf.}. It is a standard result that it satisfies (i) and (ii) \cite[29--32]{Arnold1968}.

Consider (see paragraphs 8--9) the dynamical system $(\psi,\Sigma_{\psi},\nu,T)$ where $\psi=[0,1)\times \Omega$, $\Sigma_{\psi}$ denotes  the product $\sigma$-algebra, $\nu$ is the product measure and
\begin{equation}
T(u,w)=(u+\gamma\,(\textnormal{mod}\,1),\,R(w)).
\end{equation}
 Let $V$ be the Koopman operator of $(\psi,\Sigma_{\psi},\nu,T)$. Let $d_{k,m}(w)$, $k,m\in\field{Z}$, be the orthonormal set of functions with $X(d_{k,m})=d_{k+1,m}$ and which, together with $1$, form an orthonormal basis of $L_{2}(\Omega)$ (because of (ii) above, such a set exists).
Define
$$
\bar{h}_{k}(u,w)=e^{2\pi i k u}\qquad \mbox{and} \qquad p_{l,k,m}(u,w)=e^{2\pi i l u}d_{k,m}, \quad k,l,m\in\field{Z}.
$$
$\bar{h}_{k}$ and $p_{l,k,m}$ form an orthonormal basis of $L_{2}(\psi)$ because
$d_{k,m}$ and $1$ form an orthonormal basis of $L_{2}(\Omega)$ and $e^{2\pi i k u}$ form an orthonormal basis of $L_{2}([0,1))$, and for the spaces of square integrable functions $L_{2}(Y)$ and $L_{2}(Z)$ with orthonormal bases $y_{i}(u)$ and $z_{j}(w)$, $y_{i}(u)z_{j}(w)$ form an orthonormal basis of $L_{2}(Y\times Z)$ \cite[47]{AlickiFannes}. Clearly, $V(\bar{h}_{k})=e^{2\pi i k\gamma}\bar{h}_{k}$. Thus $\bar{h}_{k}$ are proper functions with proper values $e^{2\pi i k\gamma}$, and $V$ restricted to the closed subspace $F_{0}$ generated by the $\bar{h}_{k}$ has pure point spectrum (this is point I) referred to by von Neumann). One easily verifies that  $V(p_{l,k,m})=e^{2\pi i l\gamma}p_{l,k+1,m}$. We get rid of the constant factors by setting $t_{l,k,m}=e^{2\pi i lk\gamma}p_{l,k,m}$. Clearly, $V(t_{l,k,m})=t_{l,k+1,m}$.  The set of pairs $(k,m)$ can be put into one-to-one correspondence with $\field{Z}$. If $\bar{h}_{i,j}$ is used to denote $t_{l,k,m}$ whenever $(l,m)$ corresponds to $j$, $V(\bar{h}_{i,j})=\bar{h}_{i+1,j}$. Consequently, $V$ acting on
$F_{0}^{\perp}$ has an infinitely multiple homogeneous Lebesgue measure spectrum (this is point II) referred to by von Neumann).

Recall that $h_{k}$ and $h_{i,j}$ with $U(h_{i,j})=h_{i+1,j}$ form an orthonormal basis of $L_{2}(\phi)$ and that $\bar{h}_{k}$ and $\bar{h}_{i,j}$ with $V(\bar{h}_{i,j})=\bar{h}_{i+1,j}$ form an orthonormal basis of $L_{2}(\psi)$.
Clearly, the operator defined by \begin{equation}\label{SI}
W(\sum_{k}\alpha_{k}h_{k}+\sum_{i,j}\beta_{i,j}h_{i,j}):=\sum_{k}\alpha_{k}\bar{h}_{k}+\sum_{i,j}\beta_{i,j}\bar{h}_{i,j},\,\,\, \alpha_{i},\,\beta_{i,j}\in\field{C},
\end{equation} is unitary. Hence $(\phi,\Sigma_{\phi},\mu,S)$ and $(\psi,\Sigma_{\psi},\nu,T)$ are spectrally isomorphic.

Now (see paragraph 9) we already know that for $(\phi, \Sigma_{\phi},\mu,S)$ it holds that $1''\neq 1'''$. It remains to show that for $(\psi, \Sigma_{\psi},\nu,T)$ we have $1''=1'''$, implying that $(\phi, \Sigma_{\phi},\mu,S)$ and $(\psi, \Sigma_{\psi},\nu,T)$ are \emph{not} spacially isomorphic. For $(\psi,\Sigma_{\psi},\nu,T)$ the set $1''$ consists of all proper functions of absolute value $1$. Because spectrally isomorphic systems have the same proper values, the proper functions are $\delta e^{2\pi i k u}$, $|\delta|=1$. If $g\in 1'''$, then $g(u+\gamma,R(w))=\delta e^{2\pi i k u}g(u,w)$ for some $k\in\field{Z}$. Let us expand
$g(u,w)$ in terms of the $d_{k,m}$:
$$
g(u,w)=n_{0}(u)+\sum_{k,m}n_{k,m}(u)d_{k,m}(w).
$$
Then
\begin{eqnarray}\label{ti}
g(u+\gamma,Rw)=n_{0}(u+\gamma)+\sum_{k,m}n_{k-1,m}(u+\gamma)d_{k,m}(w)\,\,\,\,\textnormal{and}\nonumber\\
\delta e^{2\pi i k u}g(u,w)=\delta e^{2\pi i ku}n_{0}(u)+\sum_{k,m}\delta e^{2\pi i k u}n_{k,m}(u)d_{k,m}(w).
\end{eqnarray}
Consequently,
\begin{equation}\label{aux:1}
n_{k-1,m}(u+\gamma)=\delta e^{2\pi i ku}n_{k,m}(u).
\end{equation}
Taking the norms of both sides of (\ref{aux:1})  yields $\left\|n_{k-1,m}\right\|=\left\|n_{k,m}\right\|$.
Recall that if $\kappa,\kappa_{j}$ are elements in an Hilbert space with $\kappa_{j}$ pairwise orthogonal, then $\sum_{j}\kappa_{j}=\kappa$ if, and only if, $\sum_{j}\left\|\kappa_{j}\right\|^{2}$ converges. Hence  $\sum_{k,m}\left\|n_{k,m}\right\|^{2}$ converges and because $\left\|n_{k-1,m}\right\| = \left\|n_{k,m}\right\|$, $n_{k,m}=0$ for all $k,m\in\field{Z}$ and $g(u,w)=n_{0}(u)$.
Equations (\ref{ti}) imply that
$$
n_{0}(u+\gamma)=\delta e^{2\pi i k u}n_{0}(u).
$$
It is well known that $\{\vartheta_{l}(u):=e^{2\pi i l u}\ :\ l\in\field{Z}\}$, form an orthonormal basis of $L_{2}((0,1])$. Hence $n_{0}(u)=\sum_{l}o_{l}\vartheta_{l}(u)$, $o_{l}\in\field{C}$. Then
\begin{equation}
n_{0}(u+\gamma)=\sum_{l}o_{l}e^{2\pi i l\gamma}\vartheta_{l}(u)\,\,\,\,\textnormal{and}\,\,\,\,\delta e^{2\pi i k u}n_{0}(u)=\sum_{l}\delta o_{l-k}\vartheta_{l}(u).
\end{equation} Consequently, $|o_{l-k}|=|o_{l}|$ for all $l$. This implies that $k=0$ because $\sum_{l}|o_{l}|^{2}$ converges and the $o_{l}$ cannot be all $0$. But for $k=0$,  $g$ is a proper function of $V$. Thus $1''=1'''$.

\section{Von Neumann's 1941 Letter: Literature Referring to the Result}\label{Literature}
To summarise, \emph{the proof in von Neumann's 1941 letter shows that for dynamical systems with mixed spectrum spectral isomorphism does not imply spacial isomorphism}.
The letter is the first document with a proof that spectral isomorphism does not imply spacial isomorphism. One would have expected that von Neumann went on to publish this result. However, this never happened, and we do not know why. In order to place the letter in historical context and to later assess in Section~\ref{Conclusion} the importance of it, we need to discuss the extant mathematical literature that refers to the result of the letter.

Two of Halmos' publications are relevant here. First,  \citeasnoun{Halmos1949}, in a survey of recent advances in ergodic theory, remarks that for mixed spectrum systems it can be shown that spectral isomorphism does not imply spacial isomorphism (the remark is three paragraphs long). To quote: ``This construction has not been published so far---it is the result of joint work by von Neumann and myself'' \cite[1025]{Halmos1949}. No proof is presented but a few details are given, which make clear that Halmos refers to the result of the letter. Namely, he states that the dynamical systems which are spectrally but not spacially isomorphic are the systems $(\phi,\Sigma_{\phi},\mu,S)$ and $(\psi,\Sigma_{\psi},\nu,T)$ mentioned in the letter (see Section~\ref{Explanation}). Moreover, he states that the proof is based on the newly-introduced notions of generalised proper values and generalised proper functions. Halmos' (1949) remark is important for three reasons: it gives a few details of the proof, it shows that the result in the letter was joint work by von Neumann and Halmos and that they intended to publish it, and it suggests that generalised proper values and generalised proper functions were first introduced by von Neumann and Halmos.

Second, the only publication we have found where a similar proof of the result of the letter is given is relatively late in \citeyear*{Halmos1956} in Halmos' `Lectures on Ergodic Theory' (in the Chapter `Generalized Proper Values'). There are no remarks to the effect that a similar construction was used by von Neumann and Halmos to first show that spectral isomorphism does not imply spacial isomorphism. (Halmos might not have felt the need to include such commentary in lecture notes.)
As already stated, Halmos' proof is similar to the one in the letter, but there are also differences, most importantly the following: first, the requirements on $R(w)$ differ. Von Neumann requires that $(\Omega,\Sigma_{\Omega},\rho,R)$ has (i) pure continuous spectrum and that (ii) the Koopman operator acting on $H_1^{\perp}$, where  $H_1$  is the subspace generated by the function 1, has infinitely multiple homogenous Lebesgue measure spectrum. The example he gives is the two-shift on the symbols 1 and -1. \citeasnoun{Halmos1956} requires $R(w)$ to be more specific and different from the two-shift, namely to be a mixing system on the unit circle such that (i) and (ii) hold.\footnote{Clearly, this difference is inconsequential because all that is needed for the proof are (i) and (ii).} Second, von Neumann calculates that $1'''\neq 1''$ for $S(u,v)$ and that $1'''=1''$ for $T(u,w)$. Halmos instead focuses on the least positive integer, not ruling out infinity, such that $1$ with $n$ dashes equals $1$ with $n+1$ dashes, and he calculates that this number is $3$ for $S(u,v)$ and $2$ for $T(u,w)$. Third, Halmos explicitly constructs the spectral isomorphism. Von Neumann only shows that the spectra of $(\phi,\Sigma_{\phi},\mu,S)$ and $(\psi,\Sigma_{\psi},\nu,T)$ are the same, and he seems to have inferred from this that they are spectrally isomorphic, which is easy to show in this case. (Yet in general it is not true that systems with the same spectrum are spectrally isomorphic -- see \citename{Halmos1951} 1951, 75.) Finally, unsurprisingly, Halmos's proof is relatively detailed whereas von Neumann's proof in the letter is sketchy.

Another publication referring to the result of the letter is \citeasnoun{Anzai}. Among other things, this paper proves the same result as the letter, viz.\ that spectrally isomorphic dynamical systems with mixed spectrum need not be spacially isomorphic. However, both the non-isomorphic dynamical systems and the methods used to prove the result are entirely different from the ones in letter. Interestingly, in the introduction one finds the following acknowledgements:
\begin{quote}
The author is much indebted to Professor S. Kakutani for his kind
discussions [...]. Further he taught the author that Professor J. von Neumann had proved
the following theorem: The ergodic transformation $T(x, y)\rightarrow (x+\gamma, x+y)$
on the torus is spectrally isomorphic to the direct product transformation of
the translation $x \rightarrow x+\gamma$ on the circle and the shift-transformation on the infinite dimensional torus 2, though these transformations are not spacially
isomorphic to each other. This fact has been the stimulation in obtaining
the results of \S6. \cite[84]{Anzai}
\end{quote}
Clearly, Anzai is here referring to the result of the letter. His quote makes clear that
von Neumann had told Kakutani about that result. Von Neumann's letter contains a paragraph about Kakutani (see the paragraph beginning with `Re Kakutani'), telling us that Kakutani has decided to stay in the US and that if the war continues, they will try to help him as much as they can.

\citeasnoun{AbramovRussian}, a paper in Russian published in English two years later \cite{Abramov}, contains a proof of the result that totally ergodic\footnote{A dynamical system $(\phi,\Sigma_{\phi},\mu,S)$ is totally ergodic if $(\phi,\Sigma_{\phi},\mu,S^{t})$ is ergodic for all $t\in\field{Z}\setminus\{0\}$.} dynamical systems are spacially isomorphic if, and only if, their generalised proper values and generalised proper functions are equivalent (see Section~\ref{Kolm} for more on this). In the introduction he refers to \citeasnoun{Halmos1949} and remarks that generalised proper values and generalised proper functions ``were introduced by von Neumann and Halmos [...], who proved with these concepts the existence of spectrally equivalent but metrically nonisomorphic automorphisms with mixed spectrum'' \cite[37]{Abramov}. Clearly, he refers here to the result in the letter but no further details about the proof are given and there is no indication that he knew more than what was announced in \citeasnoun{Halmos1949}.

Finally, \citeasnoun{Rohlin1960} briefly mentions generalised proper values and generalised proper functions and states that in this way one can show that spectrally isomorphic systems need not be spacially isomorphic. However, no further details are given and there is no indication who proved the result. Because generalised proper values and functions are mentioned, it seems likely that Rohlin refers to the result of the letter, but it remains unclear what he knew about it.

Let us now report the main results on the isomorphism problem after the publication of the letter.

\section{The Main Results After Von Neumann's 1941 Letter}\label{Kolm}
From 1942 until today a large number of various papers on the spacial isomorphism problem have been published, and much research is still being done on this problem nowadays. Thus we can only provide a summary of the most important results on the isomorphism problem after the publication of the letter. Our discussion will be chronological, and we will particularly focus on those results that are relevant to assess the importance of the letter later in Section~\ref{Conclusion}.

In 1942, shortly after the letter had been written, von Neumann and Halmos published a paper with a further contribution to the isomorphism problem. Their main result is about the same class of systems as Theorem 1, viz.\ ergodic systems with pure point spectrum.
\begin{theorem}\label{T2}
Every ergodic dynamical system $(\phi,\Sigma_{\phi},\mu,S)$ with pure point spectrum is spacially isomorphic to a rotation on a compact Abelian group \cite{NeumannHalmos}.\footnote{For  an accessible  proof,  see  \citeasnoun[46--50]{Halmos1957}. The more historically oriented papers \citeasnoun{Halmos1949}, \citeasnoun{Halmos1957}, \citeasnoun[197--198]{Mackey1974}  and Weiss (1972, 672--673) also discuss this result.}
\end{theorem}
The importance of this theorem is that it provides a normal form for the class of ergodic systems with pure point spectrum, and thus it can be used to answer many questions about this class. An interesting corollary should also be mentioned: \begin{corollary}\label{C1}
Every ergodic dynamical system $(\phi,\Sigma_{\phi},\mu,S)$ with pure point spectrum is spacially isomorphic to its inverse, i.e.\ to $(\phi,\Sigma_{\phi},\mu,S^{-1})$ \cite{NeumannHalmos}.
\end{corollary}

The next important contribution to the isomorphism problem was \citeasnoun{Anzai}.
As already mentioned in Section~\ref{Literature}, \citeasnoun{Anzai} proves the same result as von Neumann's letter, viz.\ that there are spectrally isomorphic dynamical systems with mixed spectrum which are not spacially isomorphic. However, the non-isomorphic dynamical systems as well as the construction of the proof are very different from the ones in the letter.

Up to now none of the results were on systems with pure continuous spectrum. Systems with pure continuous spectrum are generic among all dynamical systems (comeagre in the strong neighbourhood topology) \cite{Halmos1944}. Therefore, the outstanding open problem concerning the isomorphism problem was to classify systems with pure continuous spectrum, and, in particular, to find out whether for these systems spectral isomorphism implies spacial isomorphism. Kolmogorov eventually made progress on this question. According to \citeasnoun{Sinai1989}, \emph{in 1957 in a seminar in Russia Kolmogorov first presented examples of systems with pure continuous spectrum which are spectrally isomorphic but not spacially isomorphic.} This construction was then later published \cite{Kolmogorov1958,Kolmogorov1986,Sinai1959}. Let us outline Kolmogorov's proof.\footnote{This proof is standardly given to outline Kolmogorov's contribution (e.g., \citename{Halmos1961}~\citeyear*{Halmos1961}; \citename{Rohlin1967}~\citeyear*{Rohlin1967}; \citename{Weiss1972}~\citeyear*{Weiss1972}), but it is based on a definition of entropy introduced by \citeasnoun{Sinai1959}. According to \citeasnoun{Sinai1989}, this proof corresponds closely to the one Kolmogorov gave in this seminar. The proof Kolmogorov published -- \citeasnoun{Kolmogorov1958} -- was different from the one in the seminar. Most importantly, its definition of entropy was based on a theorem that turned out to be wrong \cite{Kolmogorov1986}, and so later \possessivecite{Sinai1959} definition was adopted.}

As a first ingredient, motivated by his work on information theory, Kolmogorov introduced a new spacial invariant -- nowadays called the Kolmogorov-Sinai entropy. A \emph{partition} $\alpha=\{\alpha_{i}\,|\,i=1,\ldots,n\}$ of $(\phi,\Sigma_{\phi},\mu)$ is a collection of non-empty, non-intersecting measurable sets that cover $\phi$: $\alpha_{i}\cap\alpha_{j}=\emptyset$ for all $i\neq j$ and $\phi=\bigcup_{i=1}^{n}\alpha_{i}$.
Dynamical systems and information theory can be connected as follows: each $x\in\phi$ produces, relative to a partition $\alpha$ (a coding), an infinite string of symbols $\ldots x_{-2}x_{-1}x_{0}x_{1}x_{2}\ldots$ in an alphabet of $k$ letters via the coding $x_{j}=\alpha_{i}$ if, and only if, $S^{j}(x)\in\alpha_{i},\,\,j\in\field{Z}$. Interpreting the system $(\phi,\Sigma_{\phi},\mu,S)$ as the source, the output of the source are these strings $\ldots x_{-2}x_{-1}x_{0}x_{1}x_{2}\ldots$.
If the measure is interpreted as time-independent probability, then $H(\alpha,T):=$\begin{equation}
\lim_{n\rightarrow\infty}1/n\!\!\!\!\!\!\!\!\!\!\!\!\!\sum_{i_{j}\in\{1,\ldots,k\},0\leq j\leq n-1}\!\!\!\!\!\!\!\!\!\!\!\!\!\!\!\!\!-\mu(\alpha_{i_{0}}\cap S\alpha_{i_{1}}\ldots\cap S^{n-1}\alpha_{i_{n-1}})\log(\mu(\alpha_{i_{0}}\ldots\cap S^{n-1}\alpha_{i_{n-1}}))\end{equation}
measures the average information which the system produces per step relative to $\alpha$ as time goes to infinity (\citename{FriggWerndl}~\citeyear*{FriggWerndl}; \citename{Petersen1983} 1983, pp.~233--240; \citename{Werndl2009}~\citeyear*{Werndl2009}). Now
\begin{equation}
h(T):=\sup_{\alpha}\{H(\alpha,S)\}
\end{equation}
is the \emph{Kolmogorov-Sinai entropy} of $(\phi,\Sigma_{\phi},\mu,S)$. It measures the highest average amount of information that the system can produce per step relative to a coding. It is easy to see that spacially isomorphic systems have the same Kolmogorov-Sinai entropy.

The second ingredient in the proof were Bernoulli shifts. In the modern framework of probability theory an independent process, i.e.~a doubly-infinite sequence of independent rolls of an $n$-sided die where the probability of obtaining $k$ is $p_{k}$, $k\!\in\!\bar{N}:=\{N_{1},\ldots,N_{n}\}$, with $\sum_{k=1}^{n}p_{k}\!=\!1$, is modeled as follows. Let $\Omega$ be the set of all bi-infinite sequences $w=(\ldots w_{-1},w_{0},w_{1}\ldots)$ with $w_{i}\in\bar{N}$, corresponding to the possible outcomes of an infinite sequence of independent trials. Let $\Sigma_{\Omega}$ be the set of all sets of infinite sequences to which probabilities can be assigned, and let $\mu$ be the probability measure on $\Sigma_{\Omega}$.\footnote{In detail: $\Sigma_{\Omega}$ is the $\sigma$-algebra generated by the cylinder-sets\begin{displaymath}
C^{k_{1}...k_{m}}_{i_{1}...i_{m}}\!=\!\{w\in \Omega\,|\,w_{i_{1}}\!=\!k_{1},...,w_{i_{m}}\!=
\!k_{m},\,i_{j}\!\in\!\field{Z},\,i_{1}\!<...<\!i_{m},\,k_{j}\in\bar{N},\,1\!\leq j\!\leq m\}.\end{displaymath}
The sets have probability $\bar{\rho}(C^{k_{1}...k_{m}}_{i_{1}...i_{m}})=p_{k_{1}}p_{k_{2}}\ldots p_{k_{m}}$
since the outcomes are independent. $\rho$ is defined as the unique extension of $\bar{\rho}$ to a measure on $\Sigma_{\Omega}$.} Define the shift: \begin{equation}R:\Omega\rightarrow\Omega\,\,\,\,\,\,\,R((\ldots w_{-1},w_{0},w_{1}\ldots))=(\ldots w_{0},w_{1},w_{2}\ldots).\end{equation} $(\Omega,\Sigma_{\Omega},\rho,R)$ is called a \textit{Bernoulli shift with probabilities $p_{1},\ldots,p_{n}$}.

A dynamical system $(\phi,\Sigma_{\phi},\mu,S)$ is a \emph{Lebesgue system} if there exists an orthonormal basis of $L_{2}(\phi)$ formed by the function $f(x)=1$ for all $x\in\phi$, and the functions $f_{i,j}$, $i,j\in\field{Z}$, such that $U(f_{i,j})=f_{i,j+1}$ for all $i,j$ \cite[28--30]{Arnold1968}.\footnote{$(\phi,\Sigma_{\phi},\mu,S)$ is a Lebesgue system if, and only if, $\delta 1, |\delta|=1$, are the only proper functions with absolute value 1, and (ii) $S$ acting on $H_{1}^{\perp}$, where $H_{1}^{\perp}$ is the orthogonal complement of the subspace $H_{1}$ generated by function $1$, has infinitely multiple homogeneous Lebesgue measure spectrum (cf.\ the discussion of paragraph 7 of von Neumann's letter).}
Clearly, Lebesgue systems have pure continuous spectrum and it is easy to see that they are spectrally isomorphic.\footnote{See Equation~(\ref{SI}) for the construction of the spectral isomorphism.} A simple calculation shows that all Bernoulli shifts are Lebesgue systems, and hence Bernoulli shifts are spectrally isomorphic \cite[30--31]{Arnold1968}. However, it can be calculated that the Kolmogorov-Sinai entropy of the Bernoulli shift with probabilities $p_{1},\ldots,p_{n}$ is $\sum_{i=1}^{n}p_{i}\log(p_{i})$ and thus takes a continuum of different values. Consequently:\footnote{This result is also discussed by the more historically oriented papers and chapters \citeasnoun[75--77]{Halmos1961}, \citeasnoun{LoBello}, \citeasnoun[200]{Mackey1974}, \citeasnoun[239]{ReedSimon1980}, \citeasnoun[45]{Rohlin1967},  \citeasnoun[68]{Sinai1963}, \citeasnoun[834--835]{Sinai1989} and \citeasnoun[674--676]{Weiss1972}.}
\begin{theorem}
There is a continuum of dynamical systems which are spectrally isomorphic but not spacially isomorphic.
\end{theorem}
This result by Kolmogorov was hailed as the major breakthrough since von Neumann's earlier work and stimulated much further reserach (cf.\ \citename{Halmos1961} 1961, 75--77; \citename{Rohlin1967} 1967, 45, \citename{Sinai1989} 1989, 834--835; \citename{Weiss1972} 1972, 674--676).

Another contribution which should be mentioned is \possessivecite{AbramovRussian} paper in Russian which was published in English two years later \cite{Abramov}. This paper is important because it uses generalised proper values and generalised proper functions as introduced in von Neumann's letter (cf.\ Section~\ref{Explanation}) to provide sufficient conditions for spacial isomorphism. Let us explain. For a dynamical system $(\phi,\Sigma_{\phi},\mu,S)$
denote the sequence of generalised proper functions $1''$, $1'''$, $1''''\ldots$ by $G_{\phi}^{1}, G_{\phi}^{2}, G_{\phi}^{3}\ldots$ and the sequence of generalised proper values $1'$, $1''$, $1'''\ldots$ by $H_{\phi}^{1}, H_{\phi}^{2}, H_{\phi}^{3}\ldots$.
It is easy to see that the $G_{\phi}^{n}$ and $H_{\phi}^{n}$, $n\geq 1$, are groups under multiplication \cite[57]{Halmos1956}. Generalising the notion of a pure point spectrum, \citeasnoun{AbramovRussian} introduced the definition that a dynamical system $(\phi,\Sigma_{\phi},\mu,S)$ has \emph{quasi-discrete spectrum} if $G_{\phi}=\cup_{n\geq 1} G_{\phi}^{n}$ forms a basis of $L_{2}(\phi)$.

As a next step, Abramov formalised the idea that dynamical systems have equivalent generalised proper values and generalised proper functions. Before we state the definition, note that the function $Q_{\phi}:H_{\phi}^{n+1}\rightarrow H_{\phi}^{n}$, $Q_{\phi}(f)= f(\phi(x))/f(x)$ is a group homomorphism.
Now the generalised proper values and generalised proper functions of the dynamical systems $(\phi,\Sigma_{\phi},\mu,S)$ and $(\psi,\Sigma_{\psi},\nu,T)$ are \emph{equivalent} if there exists a group isomorphism $L$ of the group $H_{\phi}=\cup_{n\geq 1} H_{\phi}^{n}$ to the group $H_{\psi}=\cup_{n\geq 1} H_{\psi}^{n}$ such that
\begin{eqnarray}
L(f)&=&f\,\,\,\textnormal{for all}\,\,\,f\in H_{\phi}^{1},\\
L(H_{\phi}^{n})&=&H_{\psi}^{n},\qquad n\in\field{N},\\
Q_{\phi}&=&L^{-1}Q_{\psi}L.
\end{eqnarray}
The main result of \citeasnoun{AbramovRussian} can now be stated.
\begin{theorem}\label{T3}
Totally ergodic dynamical systems $(\phi,\Sigma_{\phi},\mu,S)$ and $(\psi,\Sigma_{\psi},\nu,T)$ with quasi-discrete spectrum are spacially isomorphic if, and only if, the generalised proper values and generalised proper functions of the dynamical systems are equivalent.
\end{theorem}
Recall \possessivecite{Neumann1932} result that ergodic dynamical systems are spacially isomorphic if, and only if, their proper values conincide (cf.\ Theorem~1). Theorem~\ref{T3} is analogous to this result in the sense that systems with equivalent generalised proper values and generalised proper functions are spacially isomorphic. One can then also prove an analogue of Theorem 2, namely that every totally ergodic dynamical system of quasi-discrete spectrum is spacially isomorphic to an affine transformation on a compact connected Abelian group \cite{Parry1971}. However, Theorem~\ref{T3} is also very different from Theorem~1: spectrally isomorphic systems do always have equivalent generalised proper values and generalised proper functions; hence for systems with quasi-discrete spectrum spectral isomorphism does not imply spacial isomorphism.\footnote{Totally ergodic systems with quasi-discrete spectrum have either pure point spectrum or a mixed spectrum consisting of the proper values and a
infinitely multiple homogenous Lebesgue measure spectrum (cf.\ Section~\ref{Explanation}). All mixed-spectrum systems with the same proper values are spectrally isomorphic, but they need not be spacially isomorphic \cite{AbramovRussian}.} The notions of generalised proper values and generalised proper functions introduced by von Neumann in the letter proved to be very fruitful, and \possessivecite{AbramovRussian} work highlights this.

Another important contribution is \citeasnoun{Choksi}. Recall that von Neumann (1932) had proven that for \emph{ergodic} systems with pure point spectrum spacial and spectral isomorphism are equivalent (Theorem 1). It was often believed, or at least hoped \citeaffixed[495]{Neumann1932}{e.g.}, that this theorem could be extended to \emph{nonergodic} dynamical systems by invoking von Neumann's ergodic decomposition theorem, saying that any dynamical system can be decomposed in ergodic parts. \citeasnoun{Choksi} dashed this hope when he showed that for non-ergodic systems with pure point spectrum spectral isomorphism does not imply spacial isomorphism.

Finally, Ornstein's work needs to be mentioned. Since Kolmogorov's groundbreaking contribution to the spacial isomorphism problem in the late 1950s, it had been an open question whether, for certain systems, having the same Kolmogorov-Sinai entropy would be a sufficient condition for spacial isomorphism. In particular, the question arose whether Bernoulli shifts with the same Kolmogorov-Sinai entropy are spacially isomorphic. This question was answered in the positive by Ornstein's landmark work in the 1970s. In particular, Ornstein proved the following celebrated result (\citename{Ornstein1970}~\citeyear*{Ornstein1970}; see also \citename{Ornstein1974}~\citeyear*{Ornstein1974}).\footnote{This results is also discussed by the more historically oriented papers \citeasnoun{LoBello}, \citeasnoun{Sinai1989} and \citeasnoun{Weiss1972}.}
\begin{theorem}
If a Bernoulli shift with probabilities $p_{1},\ldots,p_{n}$ and a Bernoulli shift with probabilities $q_{1},\ldots,q_{m}$ have the same Kolmogorov-Sinai entropy, i.e.\ $\sum_{i=1}^{n}p_{i}\log(p_{i})=\sum_{j=1}^{m}q_{j}\log(q_{j})$, then they are spacially isomorphic.
\end{theorem}
Combined with Kolmogorov's result, this means that Bernoulli shifts are spacially isomorphic if, and only if, they have the same Kolmogorov-Sinai entropy. Ornstein developed many new techniques and so his work simulated much further research. Indeed, as already mentioned above, the isomorphism problem continues to be an active and lively research area. We have now presented the major developments on the isomorphism problem since von Neumann's 1941 letter. Now we are in a position to return to von Neumann's letter and to comment further on its historical importance.

\section{Concluding Remarks}\label{Conclusion}
We will conclude the paper by highlighting what we regard as the three major reasons why the letter is historically important.  \emph{First, von Neumann's 1941 letter to Ulam is the earliest document containing a proof of the result that spectral isomorphism does not imply spacial isomorphism}. Clearly, the letter also contains the earliest proof showing that for systems with \emph{mixed spectrum} spacial and spectral isomorphism do not always go together. Without the letter, the earliest document establishing that spectral isomorphism does not imply spacial isomorphism would be \citeasnoun{Anzai}. This led some like \citeasnoun{LoBello} to claim, erroneously as we know now from von Neumann's letter, that it was \citeasnoun{Anzai} who first established this result. From our discussion it is clear that instead Halmos and von Neumann should be credited with this. Furthermore, the letter is also the earliest document in which the notions of generalised proper values and generalised proper functions are introduced. As we have seen in Section~\ref{Kolm} (see, in particular, the work by \citename{AbramovRussian}~\citeyear*{AbramovRussian}), these notions proved fruitful: they were later employed in many papers to make progress on the isomorphism problem.

\emph{Second}, as discussed above, \citeasnoun{Halmos1949} remarks that he and von Neumann have first shown that spectral isomorphism does not imply spacial isomorphism, and \citeasnoun{Anzai} and \citeasnoun{AbramovRussian} briefly remark that von Neumann and Halmos proved that for mixed spectrum systems spectral isomorphism does not imply spacial ismorphism.  However, without the letter, apart from \possessivecite{Halmos1949} comment that the proof relies on generalised proper values and proper functions, it remained unclear how the proof proceeded. To be sure, there exists one document containing a proof of the result of the letter based on generalised proper functions and proper values, viz.\ \citeasnoun{Halmos1956}, and one could have guessed that this was von Neumann and Halmos' original proof. However, there is no commentary in \citeasnoun{Halmos1956} explaining how his proof relates to the one that he and von Neumann had intended to publish. \textit{Consequently, some uncertainty remained about the construction of the original proof. Von Neumann's letter removes this uncertainty and presents us with the original proof.}

Related to this, most of the extant articles on the history of ergodic theory do not address the question of \emph{when} it was first shown and \emph{who} first showed that spectral isomorphism does not imply spacial isomorphism (e.g., \citename{Mackey1974}~\citeyear*{Mackey1974}; \citename{Mackey1990}~\citeyear*{Mackey1990}; \citename{ReedSimon1980}~\citeyear*{ReedSimon1980}, Section VII.4; \citename{Sinai1989}~\citeyear*{Sinai1989}; \citename{Weiss1972}~\citeyear*{Weiss1972}).
The presentation of the history of the isomorphism problem in some of these publications potentially leave one with the impression that it was Kolmogorov who first proved this result when he showed that systems with Lebesgue spectrum need not be spacially isomorphic (Kolmogorov's contribution was discussed in Section~\ref{Kolm}) (\citename{ReedSimon1980}~\citeyear*{ReedSimon1980}, Section VII.4; \citename{Weiss1972}~\citeyear*{Weiss1972}).
Indeed, as mentioned above, the only paper we have found in which a correct answer to the ``when and who?" question is given is the mathematics paper \citeasnoun{Halmos1949}. It seems likely that a major reason why this question has often not been addressed is the following: Halmos' (1949), \possessivecite{Anzai} and \possessivecite{AbramovRussian} remarks are very brief and left unclear how exactly the original proof proceeded, and hence this episode of the history of the ismorphism remained obscure. We hope that the publication of the letter and the story around it will make this important episode more widely known and more widely discussed.

\emph{Third}, the letter highlights von Neumann's contribution to the isomorphism problem. Up to the present day, this contribution has  usually been taken to be (i) the introduction of the notion of isomorphism and the formulation of the isomorphism problem \citeaffixed{Neumann1932}{cf.}; (ii) the results in \citeasnoun{Neumann1932}, in particular, the proof that for ergodic dynamical systems with pure point spectrum spacial and spectral isomorphism are equivalent (cf.\ Theorem~\ref{T1}); and (iii) the results in \citeasnoun{NeumannHalmos}, in particular, the proof that any ergodic pure point spectrum system is spacially isomorphic to a rotation on a compact separable Abelian group (cf.\ Theorem~\ref{T2}) (cf.\ \citename{Halmos1957}~\citeyear*{Halmos1957}; \citename{Sinai1989}~\citeyear*{Sinai1989}; \citename{Weiss1972}~\citeyear*{Weiss1972}). \emph{The letter shows that this is not all: von Neumann, together with Halmos, was also the first to prove that spectrally isomorphic dynamical systems need not be spacially isomorphic.} Furthermore, the letter also highlights that von Neumann and Halmos were the first to introduce the important notions of generalised proper values and generalised proper functions. As discussed in Section~\ref{Kolm} (see, in particular, the work by \citename{AbramovRussian}~\citeyear*{AbramovRussian}), these notions proved very fruitful in later work on the isomorphism problem.

\bigskip
\noindent
{\bf Acknowledgement}: We wish to thank Marina Whitman, the literary heir and copyright holder of all unpublished von Neumann documents, for granting us permission to publish von Neumann's letter. We also are thankful to the Archive of the American Philosophical Society for their cooperation in providing access to the Ulam Papers.

\end{document}